\documentclass[journal,twoside]{IEEEtran}

\usepackage[cmex10]{amsmath}
\usepackage{amssymb,amsthm}

\newtheoremstyle{break}%
{9pt}{9pt}
{}
{0pt}
{\itshape}{}
{\newline}
{}
\theoremstyle{break}

\newtheorem{theorem}{Theorem}

\usepackage{mathrsfs}
\usepackage{xspace}

\DeclareFontFamily{OT1}{pzc}{}
\DeclareFontShape{OT1}{pzc}{m}{it}{<-> s * [1.100] pzcmi7t}{}
\DeclareMathAlphabet{\gyre}{OT1}{pzc}{m}{it}

\DeclareMathAlphabet{\gk}{OML}{txmi}{m}{it}

\begin{document}

\setcounter{page}{325}

\title{Transfer Functions of Generalized \\ Bessel Polynomials}
\author{JOS\'E R. MART\'INEZ, \IEEEmembership{Member, IEEE}
\thanks{Manuscript received March 1, 1976; revised February 4, 1977.}%
\thanks{The author is  with  Environmental Research and Technology, Inc., Santa Barbara, CA 93103.}}

\markboth{IEEE Transactions on Circuits and Systems, vol. CAS-24, no. 6, June 1977}%
{Martinez: Transfer Functions of Bessel Polynomials}

\maketitle

\begin{abstract}
The stability and approximation properties of transfer functions of generalized Bessel polynomials (GBP) are investigated. Sufficient conditions are established for the GBP to be Hurwitz. It is shown that the Pad\'e approximants of $e^{-s}$ are related to the GBP. 
 An infinite subset of stable Pad\'e functions useful for approximating a constant time delay is defined and its approximation properties examined. The lowpass Pad\'e functions are compared with an approximating function suggested by Budak. Basic limitations of Budak's approximation are derived.
\end{abstract}

\section{Introduction}

\IEEEPARstart{S}{ome recent} work has suggested the possibility of using ratios of generalized Bessel polynomials (GBP) to approximate the ideal delay function $e^{-s}$. \cite{Johnson}  However, the question of the stability of such rational functions was left unresolved since only a necessary, but not sufficient, condition for the stability of these polynomials was given.

The purpose of this note is to discuss several useful properties of the GBP and associated rational functions. In particular, we establish a sufficient condition for the GBP to be Hurwitz. We also show that nonminimum phase rational functions of GBP yield simultaneous maximally flat approximations of a certain order of both constant time delay and magnitude; these rational functions are shown to be related to the Pad\'e approximants of $e^{-s}$. The Pad\'e functions are compared with an approximation of $e^{-s}$  proposed by Budak. \cite{Budak}



\section{Generalized Bessel Polynomials}  \label{sec:GBP}

The GBP of degree $n$ is defined by
\begin{align}
 	B_n(s, \alpha, \beta) &= \sum_{k=1}^{n} \binom{n}{k} \frac{(n + k + \alpha - 2)^{(k)}}{\beta^k} s^{n-k},	\label{eq:gbp} 
\end{align}
where $\binom{n}{k}$ is a binomial coefficient, $(q)^{(k)} \!\!=\! q (q \!- \!1) (q\! -\! 2)$ $\cdots (q- k + 1)$, $q = n + k + \alpha - 2$, is the backward factorial function of order $k$, and $\alpha$ and $\beta$ are real parameters, $\beta \ne 0$. \cite{Burchnall}, \cite{Krall}

From \eqref{eq:gbp} it can be seen that the parameter $\beta$ is a scaling factor since
\begin{align}
 	\beta^n B_n(s, \alpha, \beta) &= B_n( \beta s, \alpha, 1),	\label{eq:scalefactor} 
\end{align}
hence the properties of $B_n$ will be determined by the parameter $\alpha$. When $\alpha = \beta = 2$, the GBP reduce to the classical Bessel polynomials of circuit theory. \cite{Storch}


\section{Hurwitz Property} \label{sec:hurwitz}


A necessary condition for $B_n(s, \alpha, \beta)$ to be Hurwitz is that its coefficients be positive. This implies that $\alpha > 1\, - \,n$ and $\beta >0$, as can be seen from \eqref{eq:gbp}. This was essentially the condition given in \cite{Johnson}. It can be readily shown that the condition is not sufficient to guarantee the Hurwitz property. \cite{Martinez1}

Theorem 1 establishes a sufficient condition for $B_n(s, \alpha, \beta)$ to be Hurwitz. It is recalled that the Hurwitz character of the special case $B_n(s, 2, \beta), \,\beta > 0$, was established long ago by Storch. \cite{Storch}
\begin{theorem}
 $ B_n(s, \alpha, \beta), n > 0$, is Hurwitz for $ \alpha \geqslant 0, \beta > 0$.
\end{theorem}
\noindent{The complete proof of Theorem 1 is found in \cite{Martinez1}.}

That the conditions of Theorem 1 are not necessary is immediately seen by referring to $B_2(s, \alpha, \beta)$. Additional information on the zeros of the GBP for $ \alpha  < 0$ as well as other general properties of the zeros are found in \cite{Martinez1}. We note, however, that for purposes of approximating $e^{-s}$, the stable GBP with $ \alpha  < 0$ generally yield poor approximations and thus are of little practical interest. \cite{Martinez1}


\section{Rational Functions} \label{sec:rational}


In this section we show that the Pad\'e  approximants of $e^{-s}$  are rational functions of GBP for integral values of
the parameter $ \alpha $ and discuss several properties of this class of rational functions. The Pad\'e  functions are then compared with Budak's approximation of $e^{-s}$ \cite{Budak}.

We shall denote the Pad\'e  approximant by 
\begin{align}
 	(n, m) &= \frac{Q_{nm}(s)}{P_{nm}(s)},	\label{eq:pamn} 
\end{align}
where $Q_{nm}$ and $P_{nm}$ are polynomials of degree $m$ and $n$, respectively. For $e^{-s}$, it is known \cite{Pade} that $Q_{nm}$ and $P_{nm}$ are given by
\begin{align}
 Q_{nm} &= 	\frac{n !}{(n+m)!} \sum_{k=0}^{m} \binom{m}{k} \frac{(n + k)!}{n!}	 \,(-s)^{m-k} \label{eq:qnm} \\
 \notag \\
 P_{nm} &=	\frac{m !}{(n+m)!}  \sum_{k=0}^{n} \binom{n}{k} \frac{(m + k)!}{m!} \,  s^{n-k}.	\label{eq:pnm}
\end{align}

Comparing $B_n(s, \alpha, \beta)$  with \eqref{eq:qnm}  and \eqref{eq:pnm}, the following theorem is immediately evident.
\begin{theorem}
 The polynomials $Q_{nm}$ and $P_{nm}$ of the $(n,m)$ Pad\'e  approximant of $e^{-s}$  are GBP. The relationship is given by
\end{theorem}
\begin{align}
 Q_{nm} &= 	\frac{n !}{(n+m)!} \, B_n(-s, \delta, 1), \quad \delta = n - m + 2, \label{eq:qnmB} \\
 \notag \\
 P_{nm} &=	\frac{m !}{(n+m)!} \, B_n(s, \alpha, 1),	\quad \alpha = m - n + 2. \label{eq:pnmB} 
\end{align}

Theorem 3 establishes the Hurwitz property of $P_{nm}$. Its proof is evident upon applying Theorem 1 to \eqref{eq:pnmB}.
\begin{theorem}
 The Pad\'e  approximant $(n,m)$ of $e^{-s}$  is stable if $m \geqslant n-2$.
\end{theorem}

It is readily ascertained that the condition of Theorem 3 is not necessary. \cite{Martinez1}, \cite{Martinez2} Moreover, it can also be shown that not all the approximants are stable. \cite{Martinez1} In fact, the number of stable as well as the number of unstable
approximants is infinite.

In order to preserve boundedness, the useful range of
physically realizable approximants is defined by $n- 2 \leqslant m
\leqslant n$, and we shall examine below the properties of the lowpass approximants $(n,n-1)$ and $(n,n-2)$.

The approximation properties of the Pad\'e  functions are  given in Theorems 4 and 5; the proofs are found in \cite{Martinez1}. It is noted that the approximation of a unit amplitude or delay is said to be maximally flat of order $k$ if, apart from the unit constant term, its Taylor expansion about the origin begins with the term $\omega^{2k}$.
\begin{theorem}[\emph{Delay Approximation}]
 The $(n,m)$ Pad\'e  approximant of $e^{-s}$  yields a maximally flat delay approximation of order $m + 1$ for $m = n- 1$ and for $m =n - 2$.
\end{theorem}
\begin{theorem}[\emph{Amplitude Approximation}]
 The $(n,n-1)$ and $(n,n-2)$ Pad\'e  approximants of $e^{-s}$  both yield a maximally flat approximation of order $n$ of the ideal lowpass characteristic.
\end{theorem}

Theorems 4 and 5 imply that the $(n,m)$ approximant  satisfies $m$ phase and $n- 1$ amplitude flatness conditions at the origin for $m = n- 1$ and $m = n- 2$. This is the maximum number of simultaneous amplitude and phase conditions that may be imposed on a transfer function of degree $m + n$. \cite{Dennis}

The usual all-pole Bessel transfer functions devote all the parameters to the delay approximation and none to approximating amplitude. By contrast, the Pad\'e  approximants allocate some parameters to phase and some to amplitude approximation. This brings us to comparing the Pad\'e  approximation of $e^{-s}$  and that suggested by Budak \cite{Budak}.

In Budak's approximation we have that
\begin{align}
 	e^{-s} &= \frac{e^{-\gamma \, s}}{e^{-(\gamma - 1)\, s}}, \quad \gamma > 0,	\label{eq:bud1} 
\end{align}
and both numerator and denominator are approximated independently by all-pole Bessel transfer functions suitably scaled and not necessarily of the same degree. The end result is a transfer function of the form
\begin{align}
 	e^{-s} \approx G_{mn}^{\gamma}(s) &= K \, \frac{B_m[2(\gamma -1) \,s, 2, 1]}{B_n(2\gamma \,s, 2, 1)},	\label{eq:bud2} 
\end{align}
where $K$ is a constant defined by $K = B_n(0, 2,1) /$ $B_m(0, 2,1)$. For $\gamma = 1$ or $m = 0$, we have the usual all-pole transfer function. Some of the attributes of \eqref{eq:bud2}   have been recently discussed in \cite{Marshak}. In the case $\gamma = 1/2$ and $m = n$ we have the well-known $(n,n)$ Pad\'e  approximant.

Marshak \cite{Marshak} has shown that the expansion of the time delay function, $t_d(\omega)$, of \eqref{eq:bud2}  has the form
\begin{align}
 	t_d(\omega) &= 1 - \gamma \left[ a_1 (\gamma \, \omega)^{2n} + \cdots \right] \notag \\
	& \qquad \quad -\left( 1- \gamma \right)  \left[ b_1 \left( 1-\gamma \right)^{2m} \, \omega^{2m} + \cdots \right],	\label{eq:td} 
\end{align}
where $a_1$ and $b_1$ are constants.  Thus for $m \ne 0$, $\gamma \ne 1$, and $m \leqslant n$, the series expansion of $t_d \left( \omega \right) - 1$ begins with the term $\omega^{2m}$. Hence the delay approximation obtained with \eqref{eq:bud2}  is maximally flat of order $m$. By contrast, the delay 
approximation given by the $(n,m)$ Pad\'e  approximant  is
maximally flat of order $m + 1$. Thus $(n,n- 1)$ and $(n,n- 2)$  
will yield better delay approximations than Budak's functions in the sense that they satisfy one more flatness condition at the origin. However, the Pad\'e  functions are 
always nonminimum phase functions, whereas Budak's
approximation also includes minimum phase functions when $\gamma \geqslant 1$.

In the sequel we will show that, except in some special cases, the lowpass functions $G_{mn}^{\gamma}(s), m< n$, yield amplitude approximations that are maximally flat of unit order for arbitrary values of $\gamma$, which is no better than is obtained with all-pole transfer functions. Moreover, it will be shown that the amplitude approximation can never be better than maximally flat of order 2, regardless of the values of $m$ and $n$. This largest order can be attained only for specific values of $\gamma > 1/2$; hence minimum-phase functions can achieve it.

Using some results from \cite{Martinez1} it can be shown that the squared-magnitude function $|G_{mn}^{\gamma}(\omega)|^2$ is given by
\begin{align}
 	|G_{mn}^{\gamma}(\omega)|^2 &= \left[ \frac{(2n)!}{(2m)!}\right]^2 \, \frac{m!}{n!}	\notag \\
	\label{eq:gmag} \\
	& \times \frac{\displaystyle{\sum_{i=0}^{m} }\binom{m}{i} \frac{ \left( 2 i \right)!}{i!} (m+i)! \left[ 2 \left( \gamma - 1 \right) \omega \right]^{2 (m-i)}} {\displaystyle{\sum_{k=0}^{n}}\binom{n}{k} \frac{ \left( 2 k \right)!}{k!} (n+k)! \left( 2 \,  \gamma \, \omega \right)^{2 (n-k)}}  \notag 
\end{align}

In order for the first nonconstant term of the Taylor expansion of \eqref{eq:gmag}  to contain the factor $\omega^{2q}, q \geqslant 1$, it is necessary that the coefficients of the terms $\omega^{2p}$, $p = 1,2, \dots, q-1$, be equal in numerator and denominator. If the coefficients of $\omega^2$ are not equal, the approximation will be maximally flat of unit order. However, only certain values of $\gamma$ will cause the appropriate pair of coefficients to be equal in \eqref{eq:gmag}. But as we shall see, with a given value of $\gamma$ one can equate only a single pair of coefficients at a time.

Equating the pair of coefficients of $\omega^{2j}, j= 1,2,... ,m$, in \eqref{eq:gmag}  we find that the equation for $\gamma$ is given by  
\begin{align}
 	\left[ \frac{\gamma}{\gamma - 1} \right]^{2j} = A_j &= \frac{\left[ (2 n)! \right]^2 \left[ (n-j)! \right]^2}{ \left(  n! \right)^2 \left[ 2 \left( n - j \right) \right]! \left( 2 n - j \right) !}  \notag \\
	\label{eq:Aj} \\
	& \qquad \times \frac{(m!)^2 \left[ 2 \left( m - j \right) \right] ! \left( 2m - j \right) !}{\left[ \left( 2 m \right) ! \right]^2 \left[ \left( m - j \right) ! \right]^2}, \notag 
\end{align}
where $j = 1, 2, \dots, m$.
Solving for $\gamma$ we obtain 
\begin{align}
 	\gamma &= \frac{\left( A_j\right) ^{1/2j}}{\left( A_j\right) ^{1/2j} \pm 1},	\label{eq:gamma} 
\end{align}
which implies that $\gamma > 1/2$. Therefore,  minimum-phase  functions can have any pair of coefficients of $\omega^{2j}$ in (11)  rendered equal.

It remains to show that coefficients can be equated only one pair at a time. To show this, suppose for the moment that a single value of $\gamma$ satisfies \eqref{eq:gamma} for more than one value of $j$, say for $j$ and $j +p,p > 0$. Then we must require that
\begin{align}
 	\left[ \frac{ \gamma^2}{\left( \gamma-1 \right)^2} \right]^j  &=	\left[ \frac{ \gamma^2}{\left( \gamma-1 \right)^2} \right]^{j+p}  \label{eq:gamma2} 
\end{align}
and this can only be true when the quantity inside the brackets is unity. This implies, in turn, that $\gamma = 1/2$. But \eqref{eq:gamma}  
only yields $\gamma > 1/2$, hence no value of $\gamma$  obtained from \eqref{eq:gamma} for a given $j$ can satisfy \eqref{eq:gamma} for any other $j$.

The implication of this property of mutual exclusion is that We can eliminate from the Taylor expansion of \eqref{eq:gmag}  only the term containing $\omega^2$. In that case the first nonconstant term will contain $\omega^4$ and the approximation will be maximally flat of order 2. This is the best that can be obtained.

The value of $\gamma$  that  causes \eqref{eq:gmag}  to yield an approximation which is maximally flat of order 2 is given by
\begin{align}
 	\gamma &= \frac{\left( 2n - 1\right)  \, \pm \sqrt{\left( 2n - 1\right) \left( 2m - 1\right)}}{2(n-m)}, \qquad m < n.	\label{eq:gamma3} 
\end{align}

 The preceding discussion indicates that the lowpass minimum-phase functions $ G_{mn}^{\gamma}(s)$ will at best satisfy a single flatness condition at the origin for values of $\gamma$  given by \eqref{eq:gamma3}. This is to be contrasted with the $(n,n-1)$ and $(n,n-2)$ approximants, both of which satisfy $(n-1)$ amplitude flatness conditions at the origin. It appears, therefore, that for a given $n \geqslant 2$ and $m = n-1$ or $m=n-2$, 
the Pad\'e  approximants will yield a better amplitude approximation  than $ G_{mn}^{\gamma}(s)$.
%


\section{Examples} \label{sec:examples}

The following examples illustrate the properties of the two types of function discussed above. For illustrative purposes we will compare the Pad\'e  approximant $(3,2)$ with $ G_{23}^{\gamma}(s)$; these are shown below:
\begin{align}
 	(3,2) &= \frac{3 s^2 - 24 s + 60}{s^3 + 9 s^2 + 36s + 60},	\label{eq:P32}  \\
	\notag  \\
	G_{23}^{\gamma}(s) &= \frac{5 (\gamma -1)^2 s^2 + 15 (\gamma -1) s + 15}{\gamma^3 s^3 +6 \gamma^2 s^2 + 15 \gamma s + 15}.	\label{eq:G23}
\end{align}

The delay functions for \eqref{eq:P32}  and \eqref{eq:G23}, respectively, are given by
\begin{align}
 	&(3,2)\!\!: \quad t_d( \omega) =	\label{eq:P32td}  \\
		& \frac{17 \,\omega^8 + 592 \, \omega^6 + 12384 \, \omega^4 + 172800 \, \omega^2 + 1440000 }{\omega^{10} + 33 \, \omega^8 + 832 \, \omega^6 + 12384 \, \omega^4 + 172800 \, \omega^2 + 1440000} \notag \\
\intertext{and} 
&G_{23}^{\gamma}(s)\!\!: 	\quad t_d(\omega) =  \label{eq:G23td} \\
& \qquad \frac{ a_4 \, \omega^8 + a_3 \, \omega^6 + a_2 \, \omega^4 + a_1 \, \omega^2 + 2025}{b_5 \, \omega^{10} + b_4 \, \omega^8 + b_3 \, \omega^6 + b_2 \, \omega^4 + b_1 \, \omega^2 + 2025},	\notag 
\end{align}
where 
\begin{align*}
 		a_4 &= 3\, (\gamma -2)  (\gamma -1)^3 \, \gamma^5 \\
		a_3 &= 9\, (\gamma -1) \, \gamma^3 \left( 4 \gamma^3 - 13 \gamma^2 + 13\,\gamma - 5\right) \\
		a_2 &= 9\, \gamma   \left( 25 \gamma^4 - 79 \gamma^3 + 120 \gamma^2 - 85\,\gamma + 25\right) \\
		a_1 &= b_1 = 135 \left( 8 \gamma^2 - 10 \gamma + 5\right) \\
		b_2 &= 9  \left( 46 \gamma^4 - 130 \gamma^3 + 165 \gamma^2 - 100\,\gamma + 25\right) \\
		b_3 &= 9 \,\gamma^2 \left( 8 \gamma^4 - 24 \gamma^3 + 32 \gamma^2 - 20 \, \gamma + 5\right) \\
		b_4 &= 3\, (\gamma -1) \, \gamma^4 \left( 3 \gamma^2 -4 \,\gamma + 2\right) \\
		b_5 &= (\gamma -1)^4 \, \gamma^6 \\
\end{align*}
Since $a_1 = b_1$ for all $\gamma$, then \eqref{eq:G23td}  is maximally flat of order 2, whereas \eqref{eq:P32td}  is flat of order 3. We cannot improve the approximation by setting $a_2 = b_2$ in (19) and solving for $\gamma$  since this results in the polynomial $(\gamma-1)^5$. Hence $\gamma = 1$, which brings us back to the all-pole function. In general, one cannot improve the approximation without reverting to the all-pole case, as can be seen from \eqref{eq:td}.

The squared-magnitude functions for \eqref{eq:P32}  and \eqref{eq:G23}  are
\begin{align}
 	| (3,2) |^2 &= \frac{9\, \omega^4 + 216\, \omega^2 + 3600}{\omega^6 + 9\, \omega^4 + 216\, \omega^2 + 3600}	\label{eq:P32mag} \\
	\notag \\
	|G_{mn}^{\gamma}(\omega)|^2 & = \frac{25 \left( \gamma - 1 \right)^4 \omega^4 + 75 \left( \gamma - 1 \right)^2 \omega^2 + 225}{\gamma^6 \omega^6 + 6\, \gamma^4 \omega^4 + 45\, \gamma^2 \omega^2 + 225}	\label{eq:G23mag}
\end{align}
It is clear that \eqref{eq:P32mag} yields an approximation that is maximally flat of order 3. Using \eqref{eq:gamma3},  we find that $ \gamma = \left( 5 + \sqrt{15} \right) /2 \approx 4.436$  causes \eqref{eq:G23mag} to have flatness of order 2, the best that we can do.

\section{Concluding Remarks} 		\label{sec:remarks}

What has the inclusion of the parameter $\gamma$ done for $G_{mn}^{\gamma}(s)$, $m < n$,  when compared to the all-pole Bessel transfer function of degree $n$? As is well known, the latter produces a delay approximation that is maximally flat of order $n$, and an amplitude approximation that is flat of unit order. On the other hand, $G_{mn}^{\gamma}(s)$ yields a delay approximation of order $m$ and, for the choice of $\gamma$ given by \eqref{eq:gamma3}, an amplitude approximation of order 2. Thus $G_{mn}^{\gamma}(s)$ can exhibit an enhanced bandwidth at the expense of the delay approximation.

Comparing the Pad\'e approximants $(n, n-1)$ and $(n, n-2)$ with $G_{mn}^{\gamma}(s)$ shows that the former yield superior approximations of amplitude and delay based on the criterion of maximal flatness. Thus for this type of lowpass function in which the degree of the numerator polynomial is lower than the denominator's by at most two units, there seems to be no reason to use $G_{mn}^{\gamma}(s)$ in preference over the Pad\'e approximants. However, this is  not  the case if one requires minimum-phase functions.

Finally, we note that we showed in \cite{Martinez1} that $(n, n-2)$ yields a delay approximation that always exhibits some time lag. This must be considered in selecting one function over another if one must choose between $(n, n-2)$ and $G_{mn}^{\gamma}(s)$.







\begin{thebibliography}{99}


\bibitem{Johnson} D. E. Johnson and J. R. Johnson,  ``On circuit-theory polynomial classes,'' \emph{IEEE Trans. Circuit Theory}, vol. CT-20, pp. 603--605, Sept. 1973.

\bibitem{Budak} A. Budak, ``A maximally flat phase and controllable magnitude,'' \emph{IEEE Trans. Circuit Theory}, vol. CT-12, p. 279, June 1965.

\bibitem{Burchnall} J. L. Burchnall, ``The Bessel Polynomials,'' \emph{Canadian J. Math}, vol. 3, pp. 62--68, 1951.

\bibitem{Krall} H. Krall and O. Frink, ``A new class of orthogonal polynomials: The Bessel polynomials,'' \emph{Trans. Am. Math. Society}, vol. 65, pp. 100--115, 1949.

\bibitem{Storch} L. Storch,  ``Synthesis of constant-time-delay ladder networks using Bessel polynomials,'' \emph{Proc. IRE}, vol. 42, pp. 1666--1675, 1954.

\bibitem{Martinez1} J. R. Mart\'inez, ``Generalized Bessel polynomials and their applications to linear phase approximation,'' doctoral dissertation, New Mexico State University, Las Cruces, NM, Oct. 1970.

\bibitem{Pade} H. Pad\'e, ``Sur la representation approch\'ee d'un fonction par des fractions rationelles,'' thesis, \emph{Ann. de l'Ecole Normale}, (3), v. 9, 1892.

\bibitem{Martinez2}  J. R. Mart\'inez, ``Pad\'e approximants, generalized Bessel polynomials, and linear phase filters,'' in \emph{Proc. 14th Midwest Symposium on Circuit Theory}, Denver, CO, 1971.

\bibitem{Dennis} F. L. Dennis and D. A. Linden, ``The derivation of pole-zero patterns by derivative adjustment,'' \emph{J. Franklin Institute}, vol. 268, p. 283--293, 1959.

\bibitem{Marshak} A. H. Marshak, D. E. Johnson, and J. R. Johnson, ``A Bessel rational filter,'' \emph{IEEE Trans. on Circuits and Systems}, vol. CT-21, pp. 797--799, Nov. 1974.






\end{thebibliography}
\end{document}